
\input amstex

\documentstyle{amsppt}
\topmatter
\title  Instantons and affine algebras I: The Hilbert scheme and vertex
operators
    \endtitle
\author I. Grojnowski \endauthor
\rightheadtext{ Instantons and affine algebras }
\address  Yale University, New Haven, CT 06520
\endaddress
\abstract
{
This is the first in a series of papers which describe the action of
an affine Lie algebra with central charge $n$ on the moduli space of
$U(n)$-instantons on a four manifold $X$. This generalises work of
Nakajima, who considered the case when $X$ is an ALE space. In
particular, this describes the combinatorial complexity of the moduli
space as being precisely that of representation theory, and thus will lead to
a description of the Betti numbers of moduli space as dimensions of weight
spaces.

This Lie algebra acts on the space of conformal blocks (i\.e\., the cohomology
of a determinant
line bundle on the moduli space \cite{LMNS}) generalising the ``insertion'' and
``deletion'' operations
of conformal field theory, and indeed on any cohomology theory.

In the particular case of $U(1)$-instantons, which is essentially the subject
of this present paper,
the construction produces the basic representation after Frenkel-Kac. Then the
well known quadratic
nature of $ch_2$,
$$ch_2 = \frac{1}{2} c_1\cdot c_1 - c_2 $$
becomes precisely the formula for the eigenvalue of the degree operator, i\.e\.
the well known
quadratic behaviour of affine Lie algebras.

}
\endabstract
\NoBlackBoxes
\date June 8, 1995 \qquad  alg-geom/9506020 \enddate
\NoBlackBoxes
\endtopmatter
\document

\define\bq{\bold Q}
\define\bz{{\bold Z}}
\define\bn{\bold N}

\define\bc{{\bold C}}

\define\aaa{\Cal A}

\define\CO{\Cal O}

\define\cstar{{\bold C^*}}

\define\uqghat{U_q\widehat{\frak g}}

\define\cv{\Cal V}
\define\cq{\Cal Q}
\define\mm{\Cal M}
\define\mmb{\bar{\Cal M}}

\define\na{\Cal M_\xi}
\define\lieg{\frak g}
\define\ghat{\widehat{\lieg}}
\define\half{{\frac{1}{2}}}

\define\hlamb{H^\Lambda}
\define\part{\Cal P}
\define\ce{\Cal E}

\define\ts{\widetilde{\Sigma}}
\define\hilb#1{\widetilde{S^{#1}X}}
\define\hilbs#1#2{\widetilde{S^{#1}_{#2}X}}
\define\tv#1{\Cal T_{\cv,{#1}}}

\define\quot{\frak Q \frak u \frak o \frak t}

\define\qbinom#1#2{\thickfracwithdelims[]\thickness0#1#2}

\head
Introduction
\endhead

This is the first in a series of papers devoted to describing the action of an
affine Lie algebra on the
moduli space of instantons on an algebraic surface $X$. This paper, which is
only an announcement, is concerned
with the ``boundary'' of moduli space; the subsequent papers will describe the
action on the interior.

We describe the idea briefly. Let $X$ be an algebraic surface, $\mm$ the moduli
space of $U(c)$-instantons on $X$
(see below for precise definitions). $\mm$ is not connected; it decomposes into
$\mm = \coprod \mm_{c_1,ch_2}$,
where $\mm_{c_1,ch_2}$ denotes those instantons with fixed first Chern class
equal to $c_1\in H^2(X,\bz)$
and second Chern character equal to $ch_2 \in \bq = H^4(X,\bq)$.

Let $\Sigma\subseteq X$ be an algebraic curve. Associated to $\Sigma$ we have
various correspondences
$ \mm @<<< \Cal P_\Sigma @>>> \mm$. Such a correspondence induces  maps
$H(\mm) \to H(\mm)$ for any cohomology
theory $H$. To describe these individual maps is very complicated. However,
these maps satisfy very simple
commutation relations, namely those defining a Lie algebra. Hence, shifting our
point of view slightly,
we see that $H(\mm)$ is a representation of this Lie algebra. This explains the
complexity of the individual
maps---they are the same as the (known) complexity of describing the action of
a Lie algebra on the individual
weight spaces of a representation.

What then needs to be described is the correspondences, the Lie algebra they
generate, and which representations occur.
The Lie algebras are {\it affine} Lie algebras, defined by the lattice
$H^2(X,\bz)$, or various sublattices such as
the lattices of algebraic cycles (Neron-Severi group). These lattices become
the ``finite part'' of the weight lattice,
with the degree operator taking value in $H^4(X,\bz)$.
(The presence of {\it affine} Lie algebras, whose characters are known to be
modular forms, is reassuringly
consistent with the remarkable work \cite{VW}, which predicts this behaviour
because of $S$-duality).

The description of the representation will have to wait for a future paper. I
hope that it is irreducible, and
determined by the K\"ahler cone and Kronheimer-Mrowka basic classes (choice of
``positive roots'' and ``highest
weight''). However, the most important invariant of a representation is
its central charge, and in this paper we show that this charge is precisely the
rank of the instanton
(i\.e\. $U(n)$-instantons give rise to level $n$ representations).

The most basic example of a correspondence is the ``elementary modifications'';
i.e\. given a divisor
$ i : \Sigma \hookrightarrow X$ we modify a vector bundle along $\Sigma$, that
is consider the correspondence
$$ \Cal P_\Sigma^n = \{ 0\to \aaa_1\to\aaa_2\to i_*\Cal E \to 0\mid \aaa_r \in
\mm, \Cal E \in Pic^n\Sigma \} $$
where $Pic^n\Sigma$ is the moduli space of holomorphic line bundles on $\Sigma$
with degree $n$.

In the case where $X$ is a curve, divisors are points, and the
analogous correspondences are precisely the geometric Hecke operators
of Drinfeld. Points however cannot interact, whereas curves on surfaces most
definitely do---their interaction
being precisely described by the lattice $H^2(X,\bz)$. It was the basic
observation of Nakajima, in the case
of an ALE space $X$, that the interactions of the correspondences are described
by the Serre relations.

Variants of these correspondences $\Cal P_\Sigma^n$ are the subject of the
sequel to this paper.

This paper is concerned with a simpler correspondence, which removes a point
from an instanton to produce a new one. As a
vector bundle modified along a point becomes only a torsion free sheaf, and not
a vector bundle, one should think
of this correspondence as acting along the ``boundary'' of the moduli space. We
find that the algebra generated by
this correspondence is essentially an affine Heisenberg Lie algebra (\S3, 5).

Obviously, these correspondences are well known in the literature. For example,
the elementary modifications
along a divisor appear in \cite{MO}, where they stratify $\mm$ by what strongly
looks like paths to the
highest weight vector of a representation (i\.e\. by the ``crystal'' basis),
and notably in \cite{KM}
where they are used to impose enough relations on the Donaldson polynomials to
determine them in terms of
certain basic classes.

Even our simple correspondence of inserting a point is a common technical tool;
see for example \cite{GL}.
Thus our main contribution is to insist that one should study the algebra of
these correspondences, and that this is easy.

Finally, these correspondences act on  {\it any} cohomology theory. In the most
intersting case,
the cohomology of $\mm$ with coefficents in a determinantal line bundle, one
produces the action of an affine algebra
on the space of conformal blocks\footnote{This paragraph is the consequence of
conversations with Greg Moore; see \cite{LMNS}.}.
This space, for which there is now a dimension formula \cite{LMNS}
generalising that of Verlinde for curves, has staggering implications for
representation theory.

Acknowledgements: It should be clear that this paper is inspired by Nakajima's
fantastic work \cite{Na}.
It is a pleasure to acknowledge helpful conversations with A\. Beilinson, R\.
Dijkgraaf, L\. Fastenberg,  I\. Frenkel,
D\. Gieseker, P\. Kronheimer, Jun Li, A\. Losev,  G\. Moore, N\. Nekrasov and
S\. Shatashvili;
and support both intellectual and moral from I\. Frenkel and G\. Moore,
without which this work would not have occurred. Portions of these results were
announced at talks
at UCLA in November 1994 and UNC-Chapel Hill in April 1995.

\head
1. Algebraic Preliminaries
\endhead

Let $V$ be a complex vector space, $t^{-1}V[t^{-1}] = V
\otimes_{\bc}t^{-1}\bc[t^{-1}]$ the associated space
of loops at $V$ which vanish at $\infty$. We make this a graded vector space by
setting
$\text{deg}(v\otimes t^n)= -n $. We write $v_n $ for $v\otimes t^n$.

Let $S = S(t^{-1}V[t^{-1}])= \oplus S_n$ be the graded polynomial algebra in
infinitely many variables.
We make this a Hopf algebra by defining $\Delta v_n = v_n \otimes 1 +1 \otimes
v_n$, for $v\in V$.
Then $S$ is a free commutative and cocommutative Hopf algebra.
Conversely, given such a Hopf algebra $S$, we can reconstruct $V\otimes t^{-n}$
uniquely as the
space of primitive elements of degree $n$.

Now suppose we are given a graded symmetric bilinear form $(,)$ on $S$ such
that multiplication and
comultiplication are adjoint, i\.e\. such that $(S_n, S_m) = 0$ if $n\neq m$,
and $(xy,z)= (x\otimes y,\Delta z)$
(where $S\otimes S$ inherits a bilinear form by $(x\otimes y,a\otimes
b)=(x,a)(y,b)$).
Such a form is completely specified by its values on the primitive elements,
i\.e\. by the values
$(\alpha^i_n,\alpha^j_n)$, where ${\alpha^i}$ runs through a basis of $V$.

Given a non-degenerate such form, we can define the action of $tV[t] = V\otimes
t\bc[t]$ on $S$,
by defining $v_n$ to be the adjoint of $v_{-n}$ for $n>0$, i\.e\.
$$(v_nx,y)=(x,v_{-n}y), \qquad n>0.$$
As $\Delta$ is an algebra homomorphism, it follows that $v_n$ acts as a
derivation on $S$,
$$\multline
 (v_n(xy),z) = (xy,v_{-n}z)= (x\otimes y,\Delta(v_{-n}z))
 = (x\otimes y,(v_{-n}\otimes 1 + 1 \otimes v_{-n})\Delta z) \\ = (v_nx\otimes
y + x\otimes v_ny, \Delta z)
= (v_n(x)y + xv_n(y),z) \endmultline$$
and hence that the Heisenberg Lie algebra $(\oplus_{n\neq 0}V\otimes t^{n})
\oplus \bc$ acts on $S$, where
$[v_n,v_m] = \delta_{n,-m}. (v_n,v_n)$ if $n > 0$.

Given $v\in V$, we define new elements $h^v_n \in S_n$ by
$$ H^v(t)= \sum_{n\geq 0} h^v_nt^n = \exp\left(\sum_{n\geq 1} v_n
t^n/n\right).$$
Then if the elements $\alpha^i$ form a basis of $V$, it is well known that the
elements $h^{\alpha^i}_n$
are algebraically independent and generate $S$.
As $t^{-1}V[t^{-1}]$ consists of primitive elements, the $H^v(t)$ are
`group-like', i\.e\.
$$ \Delta h_n^v = \sum_{0\leq a \leq n} h^v_a \otimes h^v_{n-a} $$
and the inner product is given by
$$ \sum_{n,m \geq 0} (h^v_n,h^w_m) t^ns^m = \exp\left(\sum_{n\geq 1}
(v_n,w_n)/n \cdot (ts)^n/n\right). $$

This infinite family of Heisenberg Lie algebras just constructed is
still rather flabby; however inside this space of algebras
(parametrised by maps from $\bz_+$ to non-degenerate quadratic forms
on $V$) there are certain remarkable families with much larger
symmetries; namely the vertex algebras \cite{B,FLM} and $q$-vertex algebras
\cite{FJ}.

We suppose given a lattice $L$ with non-degenerate symmetric even bilinear
form, i\.e\. $(\alpha,\alpha)\in 2\bz$ for
$\alpha \in L$, and put $V= L\otimes_\bz\bc$.
Define $(,)$ on $t^{-1}V[t^{-1}]$ by
$$ (v_n,w_m) = n (v,w) \delta_{n,m}. $$
With this inner product,
$$ \sum_{n,m \geq 0} (h^v_n,h^w_m) t^ns^m = (1-ts)^{-(v,w)}. $$
We call $S$ the {\it Fock space modeled on the lattice L}.

We also suppose given a two-cocycle $\epsilon : L \times L \to \bz/2\bz$, and
define the group algebra
of $L$ twisted by $\epsilon$, $\bc\{L\}$, as in \cite{FK,FLM}.
Define $\Cal F = S \otimes_\bc \bc\{L\}$, a $\bz_+ \times L$ graded vector
space which carries an action
of both the Heisenberg Lie algebra  $(\oplus_{n\neq 0}V\otimes t^{n}) \oplus
\bc$  and $\bc\{L\}$.
Let $V$ act on $\bc\{L\}$ by $v\cdot e^\lambda = (v,\lambda) e^\lambda$, where
$\lambda \in L$, $e^\lambda$ denotes
the corresponding element of the group algebra, and $v\in V$.

Then in these circumstances we have the well known result that a far larger
algebra acts on $\Cal F$, namely

\proclaim{Theorem \cite{FLM,B}}
$\Cal F$ is a vertex algebra, and if $L$ is positive definite, a vertex
operator algebra \footnote{
If $L$ is $\bz/2$-graded we may also make all these definitions, as long as we
work in the $\bz/2$-graded
category. So $S$ is the free super-commutative algebra on $t^{-1}V[t^{-1}]$,
i\.e\. a tensor product of an
exterior algebra and a polynomial algebra, $\Cal F$ is a super-vertex algebra,
etc.}
\endproclaim
For example, in the particular case where $L$ is positive definite and spanned
by the roots
$ \Delta = \{\alpha \in L \mid (\alpha,\alpha) =2 \} $,
then the space of vectors of conformal weight $1$ is isomorphic to the simple
Lie algebra $\lieg$ with roots $\Delta$,
and $\Cal F$ is the basic representation of $\ghat$ \cite{FK}.

If $L$ is of arbitrary signature things become much more complicated.

$$ $$
At present, there is apparantly no general definition for a quantum vertex
algebra. But if $L$ is positive
definite and spanned by the roots, then we can follow \cite{FJ} and define, for
$c\in \bn$
$$ (\alpha^i_n,\alpha^j_m) = n\delta_{nm} [nc(\alpha^i,\alpha^j)]/[n] $$
where $[n] = \frac{q^n-q^{-n}}{q-q^{-1}}$, and $\alpha^1,\dots,\alpha^l$ are a
basis of simple roots.
(There is an additional choice here, that of a positive cone in $V$).
When $c=1$ this gives $\Cal F$ the structure of the basic representation for
$\uqghat$ \cite{FJ}.

This $q$ arises in our situation when there is some ``weight'' structure on the
cohomology theory, for example a
$\cstar$-action on $X$. This does occur for ALE spaces (see \cite{N, Gr1}) but
we will stick to $q=1$ for the present paper.

\head
2. Motivic algebras
\endhead

We write this section with the minimal generality needed for this paper.
Suppose $X$ and $Y$ are two smooth proper varieties, and $Z\subseteq X\times Y$
is a subvariety, i\.e\.
a correspondence between $X$ and $Y$. We write this $X @<<< Z @>>> Y$, and
$\pi_X$, $\pi_Y$ for the two
projections from $X\times Y$ to $X$ or $Y$.

Then if $H$ is any ``reasonable'' cohomology theory we obtain honest maps
$$ R : H(X) \to H(Y), \qquad \bar{R} : H(Y) \to H(X) $$
which are adjoint with respect to the natural inner product on $H(X)$ and
$H(Y)$; $(Ra,b) = (a,\bar{R}b)$.
Here we define  $R(a) = (\pi_Y)_*(\pi_X^*a \cdot [Z])$, $\bar{R}(b) =
(\pi_X)_*([Z]\cdot \pi_Y^*b)$,
and $(a,a') = \int_* a\cdot a'$, for $a,a' \in H(X)$, $b\in H(Y)$ and where for
any space $X$, $\int : X \to pt$
denotes the projection to a point, and $[Z]$ denotes the class of $Z$ in
$H(X\times Y)$.

We mention some reasonable cohomology theories:

i) Usual homology or cohomology $H^*$; topological $K$-theory,..., cobordism,
all with (say) complex coefficients.

ii) Write $\Cal F(X)$ for the ring of constructible functions from $X$ to
$\bc$. If $\pi:X\to Y$ is a map,
$f\in \Cal F(Y), g\in \Cal F(X)$, define $(\pi^*f)(x)= f(\pi(x))$ and
$(\pi_*g)(y) = \sum_{a\in\bc}a\chi(\pi^{-1}(y)\cap g^{-1}(a))$,
where $\chi$ denotes the Euler characteristic of cohomology with compact
supports, and define
$[Z]$, for $Z\subseteq X$ a subvariety, as the characteristic function of $Z$:
$[Z](x) = 1$ if $x\in  Z$, and $[Z](x) = 0$ otherwise.

iii) If $\dim Z = (\dim X + \dim Y)/2$, then
$$ H^{\half \dim X}(X) \overset R \to{\underset  {\bar{R}} \to \rightleftarrows
}  H^{\half \dim Y}(Y)$$
where $H^*$ is the usual cohomology. If $X$ is compact K\"ahler (respectively
complex algebraic or symplectic)
we can consider the subspace of $H^{\half\dim X}(X)$ spanned by the
$(p,p)$-classes
(respectively, algebraic or Lagrangian cycles). Denote any of these subspaces
$\hlamb(X)$.
As long as $Z$ is algebraic or Lagrangian as appropriate, the topological $R$,
$\bar{R}$ preserve
these subspaces and the theories $\hlamb(X)$ are ``reasonable'' cohomology
theories.

$$ $$
As part of our definition of reasonable we require that the Kunneth map
$H(X)\otimes H(Y) \to H(X\times Y)$, $a\otimes b \mapsto \pi_X^*a \cdot
\pi_Y^*b$ is an
isomorphism. Here, if $H$ is $\bz/2$-graded then we take $\otimes$ in the
$\bz/2$-graded sense also (as in example (i)).

One can continue this list of theories as one pleases \cite{JKS}. The above
theories all take values in  vector spaces,
as our goal is to produce representations of algebras, but if the
correspondences act non-trivially on the
entire motives of $X$ and $Y$ one should also  consider functors which do not
factor through cycles homologically
equivalent to zero.

An example not on this list, but which I hope to return to, is homology of $X$
with coefficients in a given sheaf.
In our case below, the sheaf should be taken to be a determinant line bundle on
the moduli space of torsion
free sheaves, so that its cohomology is the space of conformal blocks. Then the
geometric Hecke operators
(quantum group symmetries) we produce act on the representation theory of the
double loop groups in some as yet
unknown way.

\head
3. Hilbert Schemes
\endhead

We recall some well known facts about the Hilbert scheme.
Let $X$ be a smooth algebraic surface, $S^nX = X^n/S_n$ the $n$'th symmetric
power of $X$.
For $n>1$, $S^nX$ is singular. Write $\hilb{n}$ for the Hilbert scheme of $X$,
i\.e\. for the
variety parameterising closed zero dimensional subschemes of $X$ of length $n$,
and let
$\pi:\hilb{n}\to S^nX$ be the canonical morphism sending a subscheme to its
support.
Then $\hilb{n}$ actually exists as a separated variety; it is smooth of
dimension $2n$
\cite{Gro,Fo}, $\pi$ is proper and produces a desingularisation of $S^nX$.

If $X$ is symplectic then so is  $\hilb{n}$, if $X$ is hyper-K\"ahler than
$\hilb{n}$ is the
hyper-K\"ahler resolution of the stack (orbifold) $S^nX$.
We phrase everything below in terms of the variety $\hilb{n}$, but it is often
much better to work
directly with the smooth stack $S^nX$.

If $x\mapsto nx$ denotes the diagonal map $X\to X^n \to S^nX$, then
$\pi^{-1}(nx)$ is irreducible,
and of dimension $n-1$ \cite{Br}. This fact was used in \cite{GS} to compute
the Hodge numbers of $\hilb{n}$ in
terms of those of $X$; a description of the Euler numbers which is in the
spirit of this paper can
be found in \cite{VW}.

Let $\part_n$ denote the set of partitions of $n$. If $\alpha =
(1^{\alpha_1}2^{\alpha_2}\cdots) \in\part_n$,
so $\sum_i i\alpha_i = n$, write $\ell(\alpha) = \sum_i\alpha_i$.
We have an obvious stratification of $S^nX$ by $\part_n$; the strata
$S^n_\alpha X$ has complex dimension
$2\ell(\alpha)$. Write $\hilbs{n}{\alpha}$ for $\pi^{-1}(S^n_\alpha X)$. Then
$\aaa \in \hilbs{n}{\alpha}$ if
$\aaa $ is isomorphic to a direct sum $\oplus\aaa_{i,r}$, where $1\leq r \leq
\alpha_i$, each
 $\aaa_{i,r} \in \hilb{i}$ has support a single point with multiplicity $i$,
$\pi(\aaa_{i,r}) = i\gamma_{i,r}$,
and the points $\gamma_{i,r}$ are distinct. We have $\dim \hilbs{n}{\alpha} = n
+ \ell(\alpha)$.

The open strata $S^n_{(1^n)}$, $\hilbs{n}{(1^n)}$ we also denote $(S^nX)^0$,
$(\hilb{n})^0$;
$\pi$ restricted to  $(\hilb{n})^0$ is an isomorphism.

$$ $$
Define
$$ \multline
\Lambda^0 = \Lambda^0_{ab} = \{ (\aaa_1,\aaa_2,\aaa_3)\in
\hilb{a}\times\hilb{a+b}\times\hilb{b}\mid \\
\aaa_2 \in (\hilb{a+b})^0, \text{ and there is an exact sequence }
0\to\aaa_1\to\aaa_2\to\aaa_3\to 0\}
\endmultline$$
and define $\Lambda$ to be the closure of $\Lambda^0$ in
$\hilb{a}\times\hilb{a+b}\times\hilb{b}$.
Observe that

i) If $(\aaa_1,\aaa_2,\aaa_3) \in \Lambda^0$, then $\aaa_1, \aaa_2$ and
$\aaa_3$ are all in the open stratum
of their Hilbert schemes.

ii) We have $(\aaa_1,\aaa_2,\aaa_3) \in \Lambda_{ab}$ if and only if
$(\aaa_3,\aaa_2,\aaa_1) \in \Lambda_{ba}$.

iii) Writing $+ : S^aX\times S^bX \to S^{a+b}X$ for the obvious morphism,
we have $\pi(\aaa_2) = \pi(\aaa_1)+\pi(\aaa_3)$, if $(\aaa_1,\aaa_2,\aaa_3) \in
\Lambda$.
In fact $\Lambda$ is just the correspondence of varieties produced by the
``obvious'' correspondence of
stacks
$$ \Lambda^{\text{stack}} =\{ (\aaa_1,\aaa_2,\aaa_3) \mid \aaa_2 = \aaa_1 +
\aaa_3 \}. $$
This is the basic motivic object, from which everything else follows.

iv) The dimension of $\Lambda$ is $2(a+b)$, i.e. half the dimension of the
ambient space. In fact,

\proclaim{Lemma 1} If $X$ is symplectic, then $\Lambda$ is Lagrangian (where we
change the sign of
the symplectic form on $\hilb{a+b}$ in $\hilb{a}\times\hilb{a+b}\times\hilb{b}$
as is usual).
\endproclaim
Now, let $H$ be a reasonable cohomology theory as in \S2. Write
$$ S = \oplus_{n\geq 0} H(\hilb{n}). $$
Define, for $x\in \hilb{a}, y\in\hilb{b}$, the product of $x$ and $y$,
 $$xy= (\pi_{a+b})_*((\pi_a,\pi_b)^*(x\otimes y) \cdot [\Lambda]),$$
and for $z \in \hilb{a+b}$, define
$$\Delta_{ab}z = (\pi_a,\pi_b)_*([\Lambda]\cdot \pi_{a+b}^*z) \in
H(\hilb{a})\otimes H(\hilb{b})$$
and $\Delta = \sum_{a+b =n }\Delta_{ab}$.
Also define a non-degenerate inner product
$(,): H(\hilb{n}) \times H(\hilb{m})\to H(pt)$ by
$(x,y) = \delta_{nm}\int_*x\cdot y$.

\proclaim{Theorem 2}
Equipped with this multiplication and comultiplication, $S$ is a commutative
and cocommutative Hopf
algebra. In other words, multplication and comultipication are associative,
adjoint with respect to
the inner product, and (graded) commutative (here, if $H$ is $\bz/2$-graded, so
is each $S_n$).
$\Delta$ is an algebra homomorphism.
\endproclaim
The only statement that requires proof is that $\Delta$ is an algebra
homomorphism; the rest comes free
with the formalism. This is easiest proved by using the stacks $\hilb{n}$ and
the obvious correspondendence
between them; then if the cohomology theory is the ``orbifold cohomology''
(which here is
essentially $K$-theory of $X^n$, equivariant with respect to the symmetric
group) this induces our correspondence
in the homology of $\hilb{n}$. In fact, in this form the theorem makes sense
for {\it any} variety $X$ (of any dimension),
and the orbifold cohomology\footnote{i.e\. let $\Cal F = \oplus_n
K^{S_n}(X^n,\bc)$, for $X$ any variety.
Then the theorem is that $\Cal F$ is a Fock space modeled on $H^*(X,\bc)$, with
multiplication and
comultiplication defined by the obvious correspondence. Like everything else
stated here, the proof of
this fact will appear in the longer version of this paper.};
the remarkable thing about surfaces is the additional intepretation in terms of
honest seperated smooth varieties.

$$ $$
We would now like to describe generators for $S$, and identify $S$ with a Fock
space as in \S1.
$$ $$
Let $\Sigma\subseteq X$ be a curve, $S^n\Sigma$ its $n$'th symmetric power. If
$\Sigma$ is algebraic,
we can canonically identify $S^n\Sigma$ with the Hilbert scheme of zero
dimensional subschemes of $\Sigma$ of
length $n$, so we can also regard $S^n\Sigma$ as contained in $\hilb{n}$. It is
a smooth subvariety,
and Lagrangian if $X$ is symplectic. Remarkably, the following generalisation
appears to be new:

Write $\ts$ for the subspace of points $x\in \hilb{n}$ such that $\pi(x)\in
S^n\Sigma$,
and $\ts^0_\lambda$ for $\ts\cap \hilbs{n}{\lambda}$, where $\lambda\in
\part_n$, and
$\hilbs{n}{\lambda}$ is the piece of stratification defined above.
Explicitly, $\aaa \in \hilb{n}$ is in $\ts^0_\lambda$ if
$\lambda=(1^{\alpha_1}2^{\alpha_2}\cdots)$, and
$\aaa $ is isomorphic to a direct sum $\oplus\aaa_{i,r}$, where $1\leq r \leq
\alpha_i$, each
 $\aaa_{i,r} \in \hilb{i}$ has support a single point with multiplicity $i$,
$\pi(\aaa_{i,r}) = i\gamma_{i,r}$,
and the points $\gamma_{i,r}$ are distinct. Write $\ts_\lambda$ for the closure
of $\ts_\lambda^0$ in $\hilb{n}$.

\proclaim{Proposition 3} i) $\ts$ has pure dimension $n$, $\ts^0_\lambda$ has
pure dimension $n$, and
the $\ts_\lambda$, $\lambda\in \part_n$, are precisely the irreducible
components of $\ts$.

ii) If $X$ is symplectic, then $\ts$ is a Lagrangian submanifold.
\endproclaim
\comment
\demo{Proof} As the $\hilbs{n}{\lambda}$ stratify $\hilb{n}$, the
$\ts^0_\lambda$ partition $\ts$.
Each $\ts^0_\lambda$ is fibred over $S^n\Sigma \cap \hilbs{n}{\lambda}$ which
has dimension $\ell(\lambda)$.
The fibres have dimension $\sum\alpha_i(i-1)= n - \ell(\lambda)$, if
$\lambda=(1^{\alpha_1}2^{\alpha_2}\cdots)$.
Whence (i). As for (ii), observe $\Lambda$ is Lagrangian, as are curves $\Sigma
\subseteq X$; that
multiplication as defined above (via the correspondences) sends Lagrangian
subvarieties to Lagrangian
subvarieties, and that $\ts$ is just $\Sigma\cdots\Sigma$ ($n$ times). (ii)
follows.
\enddemo
\endcomment
See \S6 below for more remarks on $\ts$.

Write $h_n^\Sigma$ for the class of $S^n\Sigma \hookrightarrow \hilb{n}$ inside
$S$. (Note that as a
consequence of the proposition, we can define this class even if $\Sigma$ is
not algebraic, though we
do not need to.) We adopt the convention that $h_0^\Sigma = 1$, for all
$\Sigma$.

\proclaim{Proposition 4} i) The elements $h^\Sigma_n$ are group-like, i.e.
$$ \Delta h^\Sigma_n = \sum_{a+b=n} h^\Sigma_a\otimes h^\Sigma_b.$$
ii) Let $\Sigma$, $\Sigma'$ be two curves in $X$. Then
$$ \sum_{n,m \geq 0} (h^\Sigma_n,h^{\Sigma'}_m) t^ns^m =
(1-ts)^{-(\Sigma,\Sigma')} $$
where $(\Sigma,\Sigma')$ denotes the inner product in $H(X)$, as usual.
\endproclaim
Now, for $n \geq 0$ define $\Gamma^X = \Gamma^X_n = \hilbs{n}{(n)}$, a closed
irreducible subvariety of $\hilb{n}$.
So $\Gamma^X=\{\aaa\in\hilb{n} \mid \pi(\aaa) = nx,\text{ for some }x\in X\} =
X \times_{S^nX} \hilb{n}$.
Write $\pi:\Gamma^X\to X$, and define for any submanifold $Z\subseteq X$,
$\Gamma^Z = \Gamma^Z_n \subseteq \Gamma^X_n$
by $\Gamma^Z = \pi^{-1}(Z)$, i\.e\. so that the diagram
$$\CD \Gamma^Z @>>> \Gamma^X @>>> \hilb{n} \\
@V{\pi}VV @V{\pi}VV @V{\pi}VV \\
Z @>>> X @>{n}>> S^nX \endCD $$
is Cartesian. In particular, if $\Sigma \subseteq X$ is an algebraic curve,
$\Gamma^\Sigma = \ts_{(n)}$ in the notation
above.
Also write $r^Z_n$ for the class of $\Gamma^Z_n$ in $H(\hilb{n})$. We then have

\proclaim{Lemma 5} There exists a function $f : \bz\to\bz$, and polynomials in
infinitely many variables
$p^n_a(x_1,x_2,\dots), q^n_b(x_1,x_2,\dots)$ such that for any submanifolds
$Z,Z'$ of $X$,
$$\align
&  \text{i) }  (r^Z_n,r^{Z'}_n) = (Z,Z')f(n), \text{ and }  \\
& \text{ii) }  \Delta r^Z_n = \sum_{a+b=n} p^n_a(r^Z_1,r^Z_2,\dots) \otimes
q^n_b(r^Z_1,r^Z_2,\dots)
\endalign $$
\endproclaim

Now observe that for a fixed algebraic curve $\Sigma$, the algbera generated by
the $r^\Sigma_n$, $n\in\bz_+$,
is the same as the algebra generated by the $h^\Sigma_n$, $n\in\bz_+$. As a
consequence of lemma 5, we can write
$h^\Sigma_n$ as a polynomial in the $r^\Sigma_i$ with ``universal''
coefficients; i.e. coefficients independent of $\Sigma$.
This allows us to define $h^Z_n$, for any $Z\subseteq X$ by the same formulae.
It then follows from lemma 5 and proposition 4 that the $h^Z_n$ are also group
like, and satisfy the
 same inner product formulae as in proposition 4.

Take $H(X) = H^*(X)$ (usual cohomology), so $S = \oplus H^*(\hilb{n})$, and let
$Z_i$, $i=1,\dots,l$ run through
submanifolds of $X$ such that the classes $[Z_i]$ form a basis in $H^*(X,\bc)$.

\proclaim{Proposition 6} The elements $r^{Z_i}_n$ freely generate $S$ as an
algebra.
\endproclaim

Now let $X$ be projective, so that the lattice $H^*(X,\bz)/\text{torsion}$ is a
non-degenerate lattice
with respect to the form $(,)$. We have thus proved

\proclaim{Theorem 7} $S$ forms a Fock space modeled on the lattice
$H^*(X,\bz)/\text{torsion}$.
\endproclaim

If $X$ is affine, one easily shows that if we take $H(X) = H_{\half\dim
X}(X,\bc)$ to be middle dimensional
Borel-Moore homology, that the lattice $L = H_{\half\dim X}(X,\bz)$ is
non-degenerate and $S$ (for this cohomology
theory) is a Fock space modeled on $L$. This is precisely the case that occurs
for $X$ an ALE space.

For smooth projective $X$, let us agree to write $S^\Lambda =
\oplus\hlamb(\hilb{n})$, where $\hlamb$  as in \S2
denotes either toplogical, $(n,n)$, holomorphic or Lagrangian middle
dimensional cycles, and let us still write
$S = \oplus H^*(\hilb{n})$. Clearly $S^\Lambda$ is a Hopf subalgebra of $S$,
and we can describe it explicitly
using the above theorem.

For example, consider the Hodge decomposition of $H^*(\hilb{n})$; write
$S^{ab}_n=H^{a,b}(\hilb{n})$.
Then multiplication in $S$ preserves all this grading: $S^{ab}_n\cdot S^{cd}_m
\subseteq S^{a+c,b+d}_{n+m}$.
The generators $r_n^Z\in S_n$ have degree $(n-1,n-1) + \deg Z$, i\.e\. if
$[Z]\in H^{p,q}(X,\bc)$ then
$r_n^Z\in S_n^{n+p-1,n+q-1}$. This gives another proof of \cite{Got,3.1},
independent of \cite{GS}.

Observe that $S^\Lambda$ is usually not a Fock space modeled on a lattice. For
example, if the odd cohomology
of $X$ vanishes, the generating function $\dim S^{n,n}_n z^n$ is the coefficent
of $u^0$ in
$$ \prod_{n\geq 1}
\big((1-z^nu)(1-z^nu^{-1})\big)^{-h^{2,0}}(1-z^n)^{-h^{1,1}}$$
which is not the generating function of a Fock space.

\head
4. Vertex algebras and $U(1)$-instantons
\endhead
Let $X$ be a smooth projective surface. A torsion free sheaf $\ce$ is a
coherent sheaf of $\CO_X$-modules
which is torsion free as a $\CO_X$-module. If $\ce^*$ denotes the dual of
$\ce$, we have a canonical exact sequence
$ 0 \to \ce\to\ce^{**}\to \Cal Q\to 0$, where $\Cal Q$ is coherent of finite
length, and $\ce^{**}$ is locally free.
We say $\ce$ has rank $c$ if $\ce^{**}$ does.

Let $H$ be a fixed very ample divisor on $X$, and $\mm$ (resp\. $\mmb$) the
space of $H$-stable rank $c$
torsion free sheaves (resp\. $H$-semistable torsion free sheaves, modulo the
usual equivalence relation \cite{Gi}).
Then $\mmb$ is a projective variety, and $\mm \subseteq \mmb$ an open
subvariety \cite{Gi}.
Clearly $\mmb = \coprod_{c_1,k} \mmb_{c_1,k}$, where $\mmb_{c_1,k}$ consists of
torsion free sheaves $\ce$
with Chern classes $c_1(\ce) = c_1, ch_2(\ce)=k$.
Also, $c_1(\ce) $ lands in the lattice of algebraic cycles $c_1(\ce) \in
\hlamb(X,\bz)$.

In the sequel we will be concerned with elementary modifications, and hence the
interior of moduli space. For now,
as we are only concerned with the boundary of moduli space, let us suppose
$c=1$.
Further, for simplicity, suppose $\pi_1(X)=0$. Then $H^1(X,\bz)=H^3(X,\bz)=0$,
and $H^2(X,\bz)$ is torsion free.
Suppose also that $X$ is spin, so $H^2(X,\bz)$ is even.

Now, in this case line bundles have no moduli, and so the space of torsion free
sheaves is isomorphic to
$$ \oplus_{\lambda\in \hlamb (X,\bz)} (\oplus_n \hilb{n}) \otimes e^\lambda.$$
Observe that the moduli of torsion free sheaves $\ce$ such that $\ce^{**}$ is
isomorphic to $\Cal L$
is canonically isomorphic to those with  $\ce^{**}$  isomorphic to $\Cal L'$;
the map is just
$\ce \mapsto \ce\otimes (\Cal L'\otimes \Cal L^{-1})$. This defines the action
of the lattice $\hlamb(X,\bz)$
on the moduli space of rank 1 torsion free sheaves. (If $X$ is not simply
connected, we must replace this
action by elementary modifications; this will be explained in the sequel).

Let us write $L = H^*(X,\bz)$, and $ \mm' = (\oplus_n \hilb{n})\times L$,
$$\Cal F = \oplus H^*(\hilb{n},\bc)\otimes \bc\{L\} = S \otimes \bc\{L\}.$$
Then $\mm'$ is ``almost'' the moduli space of topological $U(1)$-instantons;
i.e\. the hyper-K\"ahler resolution
of the ideal ASD-connections on toplogical $U(1)$-bundles. (It would be
interesting to give a purely algebraic
construction of this space). It differs from this by the term $H^0(X,\bz) +
H^4(X,\bz)$,
which I do not know how to intepret geometrically. Yet.

In any case, the results of \S3 tell us
\proclaim{Theorem} $\Cal F$ is a vertex algebra.
\endproclaim
\noindent
and that if $X$ is an ALE space, so $L=H^2(X,\bz)$ is negative definite, that
the analogously defined $\Cal F$
(which is now precisely the homology of rank 1 torsion free sheaves) is just
the basic representation
of $\ghat$, where $\ghat$ is the affine Lie algebra associated to $L$.
\indent

Also, we remark that the theorem is true for an arbitrary compact K\"ahler
4-manifold, where we use the orbifold cohomology.

\head
5. Torsion sheaves of rank $c$ produce central charge $c$
\endhead

Let $\cv$ be a vector bundle of rank $c$ on $X$, and let $\tv{n}$ consist of
the stack of torsion free sheaves
 $\ce$ such that $\ce^{**}$ is isomorphic to $\cv$. Then if $\cv$ is stable
(resp. semistable), so is
any $\ce\in\tv{n}$ (though not conversely), and in that case $\tv{n}$ is a
smooth separated variety.

Define $\cq_\cv = \oplus_n H(\tv{n})$, where $H$ is any reasonable cohomology
theory. One may handle $\cq_\cv$
as one handles $S$ in \S3, by defining an action of $S$ on $\cq_\cv$, $S\otimes
\cq_\cv @>m_\cv>> \cq_\cv$
induced by the correspondence
$$ \multline \Lambda^0 =\{ (\aaa_1,\aaa_2,\aaa_3) \in \tv{a}\times\tv{a+b}
\times \hilb{b} \mid \\
\aaa_2 \in \tv{a+b}^0, \text{ and there is an exact sequence }
0\to\aaa_1\to\aaa_2\to\aaa_3\to 0\}
\endmultline$$
and $\Lambda $ is the closure of $\Lambda^0$. Here, $\tv{n}^0$ consists of
those torsion free sheaves $\ce$
which fit into an exact sequence $\ce\to\cv\to \Cal Q$, where $\Cal Q \in
(\hilb{n})^0$.

Then one may proceed exactly as in \S3, and show
\proclaim{Theorem} $\cq_\cv$ is a module for the Heisenberg Lie algebra
$ (\oplus_{n\neq 0} H^*(X,\bc) \otimes t^n) \oplus \bc$ with central charge
$c$.
\endproclaim

It is pleasant to calculate $$[h^\Sigma_1,h^{\Sigma`}_{-1}] =
c(\Sigma,\Sigma')$$ directly. This follows directly
from the easy fact that if $\ce$ is any torsion free sheaf of generic rank $c$,
and $\bc_0$ denotes
the skyscraper sheaf at a point $0\in X$, then
\proclaim{Lemma} $\dim Ext^1(\bc_0,\ce) + c = \dim Hom(\ce,\bc_0)$.
\endproclaim
In the particular case that $c=1$, this may be intepreted as
(and in fact follows from) the fact that there is one way more to add a square
to a partition than to remove
a square.
(As one may complete $X$ at 0, to get $\bc[[x,y]]$. This admits a
$\cstar$-action, such that $x^iy^j$ have distinct
weights $i+j \leq n$ . Then the fixpoints of this $\cstar$ action on  the
Hilbert scheme of length $n$ subschemes
are just in 1-1 correspondence with partitions of $n$, and this punctual
Hilbert scheme partitions into
vector bundles over these isolated fixpoints. This makes the lemma obvious in
this case, and gives yet another reason
why the Hodge theory of $H^*(\hilb{n})$ is so simple).

\head
6. Remarks on Curves
\endhead

Suppose $\Sigma \subseteq X$ is an algebraic curve. Replace $X$ with the normal
bundle to $\Sigma$ in $X$,
so that $X$ admits a contracting $\cstar$ action with fixpoints $\Sigma$, that
is if $x\in X$,
$\lim_{t\to 0} t\cdot x$ exists and is in $\Sigma$. This $\cstar$ action
induces one on $\hilb{n}$,
and we define
$$ \Cal U = \{ \aaa\in\hilb{n} \mid \lim_{t\to 0} t\cdot x \in S^n\Sigma \} $$
where $S^n\Sigma\hookrightarrow \hilb{n}$ as in \S3.

\proclaim{Proposition} i) $\Cal U$ is open in $\hilb{n}$, and $\Cal U$ is a
rank $n$ vector bundle on $S^n\Sigma$.

ii) Suppose $X= T^*\Sigma$. Then $\Cal U$ canonically identifies with
$T^*(S^n\Sigma)$.
Under this identification, the Lagrangian subvariety $\ts \cap \Cal U$
identifies with Laumon's global
nilpotent cone, a Lagrangian subvariety in $T^*(S^n\Sigma)$ \cite{La}.
\endproclaim
Because of this proposition, a perverse sheaf on $S^n\Sigma$ with nilpotent
characteristic variety (for example,
conjecturally any automorphic sheaf) gives rise to a cycle in
$\hlamb(\hilb{n})$ via the characteristic cycle map.

Its also worth remarking that if $\Sigma$ is the affine line, $S^n\Sigma$
canonically identifies with the variety of
{\it regular} conjugacy classes in $\frak g\frak l_n$ (via the characteristic
polynomial).
Thus $S^n\Sigma$, for $\Sigma$ a curve of genus $g$, which classically
\footnote{{`classically' here means after the work of Drinfeld.}}
 one regards as a genus $g$ generalisation of
regular conjugacy classes, is here generalised `microlocally' to produce two
dimensional analogues of conjugacy classes.
(As the curve $\Sigma$ varies, this
really does feel two dimensional).

Note that we do {\it not} want to consider the stack of coherent sheaves of
length $n$ here (the analogue of the
stack of all conjugacy classes in $\frak g\frak l_n$), as in our case the
central charge would vary:
$\quot_{V,r}$ has central charge the rank of $V$, where $V$ is a vector bundle
on $X$.

\head 7. Remarks on Nakajima's quiver varieties and \cite{Gr1}
\endhead
Let $\frak g$ be a Kac-Moody Lie algebra, with symmetric Cartan matrix. In
\cite{L1},
Lusztig defined a variety, the moduli space of representations of the quiver
associated to $\frak g$,
and a Lagrangian subvariety $\Lambda$ such that the middle dimensional cycles
on $\Lambda$
($\hlamb(\Lambda)$ in the notation of \S2) realises the universal Verma module
for $\frak g$.

In \cite{Na} Nakajima constructed a modified quiver variety $\na(w)$, depending
on a highest weight $w$
and a $\xi\in \frak h$, where $\frak h$ is the ``real'' Cartan subalgebra of
$\lieg$, with the
following properties:

i) If $\xi$ is generic, then $\hlamb(\na(w))$ realises the irreducible
integrable highest weight
module with highest weight $w$, and the Chevalley generators of $\lieg $ act on
$\na(w)$ by corrspondences.

ii) If $\lieg$ is of affine type, and $tr\,\xi=0$, then $\na(w)$ is the moduli
of $U(n)$-instantons on
the ALE space $X_\xi$, with monodromy at $\infty$ determined by $w$. (We refer
to \cite{KN,N} for all these
terms).

In other words, in case (i) Nakajima ``cuts down'' a Verma module to get an
irreducible highest weight module.
Unfortunately, (i) and (ii) cannot occur simultaneously\footnote{{This was
explained to me by Greg Moore.}}.
For example, $\na(0)$ is a point if $\xi$ is generic,
but if $\xi$ is generic trace free it is the Hilbert scheme on $X_\xi$.

So, if we care about the moduli of instantons on an ALE space, we must do some
extra work from \cite{Na}.
Obviously, this is the content of this paper, which complements \cite{Na} even
in the case of an ALE space.
(Using \cite{Lu2,Na} one can obtain the results above in the quiver language
directly. This will appear in \cite{Gr1}).

The point of this series of papers is to use elementary representation theory
to obtain information
about the moduli space of instantons. In the quiver variety case (for $\xi$
generic), one may reverse this, and use the geometry
of quiver varieties to obtain new information about quantum affine algebras.

Specifically, let $\lieg$ be a Kac-Moody algebra, and $\uqghat$ the associated
quantum affine algebra at central charge 0
\cite{Dr,Gr2}. If $\frak g$ is finite dimensional, we can consider the category
of finite dimensional representations
of $\uqghat$; for the definition for general $\lieg $ see \cite{Gr1}
(these representations have the property that they restrict to
a direct sum of integrable highest weight representations of
$U_q\lieg\hookrightarrow \uqghat$). Then, as was
discovered by Drinfeld, these representations are {\it not} deformations of the
analogous representations of $\ghat$
(the ``evaluation representations'' and their tensor products); smaller terms
must be added.

The reason for this is that the varieties $\na(w)$ are not zero dimensional;
i.e. $H^*(\na(w))$ is $\hlamb(\na(w))$
plus smaller terms.

Geometrically, we take as our reasonable cohomology theory
$K^{GL_W\times\cstar}(\na(w))$,
which takes full account of the geometric symmetries of $\na(w)$.
Then in \cite{Gr1} it is proved that $K^{GL_W\times\cstar}(\na(w))$, admits an
action of $\uqghat$ (see also \cite{Gr2}).
This explains the occurance of the middle homology in \cite{Na}.

Also, to continue the advertisement of \cite{Gr1}, we construct all of
$\uqghat$. Namely, we take
equivariant cohomology of the Lagrangian subvariety of $\na(w)\times\na(w)$
consisiting of pairs with the same moment map image.
This constructs a piece of $\uqghat$ (and the same variety was independantly
discovered by Nakajima in \cite{Na2}, where he
used it to construct a piece of the enveloping algebra of $\lieg$).
These pieces fit togethor via the coproduct---write $w=w'+w''$. This defines a
$\cstar$ action on $\na(w)$, with
fixpoints of the form $\na(w')\times\na(w'')$. Then one may define a coproduct
via localisation to the fixpoints, and this
coproduct fits these pieces togethor to produce $\uqghat$.

This coproduct, at the $q=1$ non-affine level, realises the map $L_{w'}\otimes
L_{w''} \to L_{w'+w''}$,
where $L_w$ is the irreducible  highest weight module for $\lieg$ with highest
weight $w$.

Finally, we describe the irreducible modules (even at roots of unity) in terms
of certain perverse sheaves with
nilpotent characteristc variety on Nakajima's moduli space $\Cal M_0(w)$.

In the case $\frak g =\frak g\frak l_n$, this construction of the algebra is
due to Ginzburg and Vasserot \cite{GV};
the moduli space is due to Beilinson-Lusztig-MacPherson \cite{BLM}; and the
coproduct  appeared in \cite{Gr3}.
That such a geometric picture of the representation theory of $\uqghat$ should
exist was conjectured by Drinfeld,
on the basis of Kazhdan and Lusztig's description of the representation theory
of affine Hecke algebras, which is similar to this.
This was explained to me by G\. Lusztig, in 1991.

\enddocument